\documentclass[11pt]{article}
\usepackage{amsfonts, amsmath, amssymb, amsthm}
\usepackage[colorlinks=true,citecolor=blue,urlcolor=blue,linkcolor=blue,bookmarksopen=false]{hyperref}

\date{}

\author{Boris Bukh\footnote{Department of Mathematical Sciences, Carnegie Mellon University, Pittsburgh, PA 15213, USA. Email: bbukh@math.cmu.edu. Supported in part by U.S. taxpayers through NSF grant DMS-1201380.}}
\title{Random algebraic construction of extremal graphs}

\newtheorem{theorem}{Theorem}
\newtheorem{lemma}[theorem]{Lemma}

\theoremstyle{remark}

\newtheorem*{remark}{Remark}
\DeclareMathOperator{\ex}{ex}                         
\DeclareMathOperator{\Binom}{Binom}                   
\newcommand*{\F}{\mathbb{F}}                          
\newcommand*{\E}{\mathbb{E}}                          
\newcommand*{\abs}[1]{\lvert #1\rvert}                
\newcommand*{\eqdef}{\stackrel{\text{\tiny{def}}}{=}} 

\begin{document}
\maketitle

\begin{abstract}
We present a motivated construction of large graphs not containing a given
complete bipartite subgraph. The key insight is that the algebraic constructions
yield very non-smooth probability distributions.
\\[1ex]\textsc{MSC classes:} 05C35, 05D99
\end{abstract}

\section*{Introduction}
A foundational problem in extremal graph theory is Tur\'an's question: how big can a
graph be if it does not contain $H$ as a subgraph? Let $\ex(n,H)$ be the maximum number of edges in any $n$-vertex $H$-free graph. 
Tur\'an himself \cite{turan_orig} determined $\ex(n,H)$ when $H$ is a clique:
He showed that if $H=K_r$, then the maximum is attained by a complete $(r-1)$-partite graph whose
parts are as equal as possible (see \cite{aigner_monthly} for six different proofs of this result). 
Erd\H{o}s--Stone and Simonovits \cite{erdos_stone,erdos_simonovits,simonovits_stability} 
showed that \emph{for every graph $H$} the largest $H$-free graph is close (in an appropriate sense) to
a complete multipartite graph on $\chi(H)-1$ parts, where $\chi(H)$ is the chromatic number of $H$.
In particular, the asymptotic formula
\begin{equation}\label{eq:ess}
  \ex(n,H)=\left(1-\frac{1}{\chi(H)-1}\right)\binom{n}{2}+o(n^2)
\end{equation}
holds. 

It is natural to interpret $1-\frac{1}{\chi(H)-1}$ as the fraction of the total number of edges in the complete
graph. When $\chi(H)\geq 3$, the fraction is positive and \eqref{eq:ess} is a satisfactory asymptotics.
On the other hand, if $H$ is bipartite, 
then $\chi(H)=2$, and the main term in \eqref{eq:ess} vanishes, leaving a notoriously hard open problem of finding an asymptotics
for $\ex(n,H)$ when $H$ is bipartite.

In this paper we focus on the best-understood class of bipartite graphs, the complete bipartite graphs. Let $K_{s,t}$
denote the complete bipartite graph with parts of size $s$ and~$t$. The following is a basic upper bound on $\ex(n,K_{s,t})$.
\begin{theorem}[Kov\'ari--S\'os--Tur\'an]
For each $s$ and $t$ there is a constant $C$ such that
$\ex(n,K_{s,t})\leq C n^{2-1/s}$.
\end{theorem}
\begin{proof}
We let $C$ be a large constant (to be specified later). 
Suppose $G=(V,E)$ is a $K_{s,t}$-free graph. 
It suffices to prove that $G$ contains a vertex of degree less than $Cn^{1-1/s}$, for then
we may remove it, and apply the induction on the number of vertices since $C(n-1)^{2-1/s}+Cn^{1-1/s}\leq Cn^{2-1/s}$.
Assume, for contradiction's sake, that $\deg(v)\geq Cn^{1-1/s}$ for all $v\in V$. 

Let $N$ denote the number of copies $K_{1,s}$ in $G$. We count $N$ in two different ways. 
On one hand, denoting by $\deg(v)$ the degree of $v\in V$, we obtain
\[
  N=\sum_{v\in V} \binom{\deg(v)}{s},
\]
the summand being the number of copies of $K_{1,s}$ with the apex~$v$. Since $\deg(v)\geq Cn^{1-1/s}$ for all $v$ and $C$
is sufficiently large in terms of $s$,
we have $\binom{\deg(v)}{s}\geq \bigl(\tfrac{1}{2}Cn^{1-1/s}\bigr)^s/s!=2^{-s}C^s n^{s-1}/s!$
and hence 
\begin{equation}\label{eq:kst_upper}
  N\geq 2^{-s}C^s n^s/s!
\end{equation}

On the other hand, if $\{u_1,\dotsc,u_s\}$ is any set of $s$ vertices, then
no more than $t-1$ vertices can be adjacent to all of these $s$ vertices, 
as $G$ is $K_{s,t}$-free. Thus
\begin{equation}\label{eq:kst_lower}
  N\leq (t-1)\binom{n}{s}.
\end{equation}

Combining \eqref{eq:kst_upper} and \eqref{eq:kst_lower} together
with the simple bound $(t-1)\binom{n}{s}\leq (t-1)n^s/s!$ yields
a contradiction unless 
$
  C\leq 2(t-1)^{1/s}.
$
\end{proof}

Despite being a sixty-year-old result with a simple proof, the K\'ovari--S\'os--Tur\'an theorem
has been improved only once, by F\"uredi \cite{furedi_bound} who improved the bound on the constant~$C$.
Is Kov\'ari--S\'os--Tur\'an bound tight? It is for $K_{2,2}$, and $K_{3,3}$ \cite{erdos_renyi_sos,brown_constr}, but
no constructions of $K_{s,s}$-free graphs with $\Omega(n^{2-1/s})$ edges are known for any $s\geq 4$.
There are however constructions \cite{kollar_ronyai_szabo,alon_ronyai_szabo,bbk} 
of $K_{s,t}$-free graphs with  $\Omega(n^{2-1/s})$ edges
when $t$ is much larger than $s$. The aim of this paper is to present a new construction that uses
both the algebra and probability. The construction is inspired by the construction in \cite{bbk}.
As a motivation, we first explain why the standard probabilistic argument is insufficient.

\section*{Sketch of a probabilistic construction}
Our graphs will be bipartite, with $n$ vertices in each part.
We refer to the parts as `left' and `right', and denote them $L$ and $R$ respectively. 

Set $p=n^{-1/s}$. For each pair of vertices $u\in L$, $v\in R$
declare $uv$ to be an edge with probability $p$, different edges being independent.
The expected number of edges in $G$ is $pn^2=n^{2-1/s}$. 
As the number of edges in $G$ is binomially distributed, $G$ will have at least
$\tfrac{1}{2}n^{2-1/s}$ edges with probability tending to $1$ as $n\to\infty$.

We shall show that $\Pr[K_{s,t}\subset G]\to 0$ for a suitably large $t$. For a set $U$ of vertices,
let 
\[
  N(U)\eqdef \{v \in V : uv \text{ is an edge for all }u\in U \}
\]
be the common neighborhood of~$U$. 

Let $U$ be any set of $s$ vertices. We shall bound $\Pr[\abs{N(U)}\geq t]$.  
If $U$ is not contained in either $L$ or $R$, then $N(U)=\emptyset$ and $\Pr[\abs{N(U)}\geq t]=0$.
Suppose $U\subset L$ (the case $U\subset R$ is symmetric). It is clear that $v\in N(U)$ with
probability $p^s=1/n$ for each $v\in R$, and these events are independent for different $v$'s.
Hence, $\abs{N(U)}$ is a binomial random variable $\Binom(n,1/n)$. It follows
that $\abs{N(U)}$ is distributed approximately as a Poisson random variable with mean $1$,
and in particular it can be shown that $\Pr\bigl[\abs{N(U)}\geq t\bigr]\leq 1/t!$.

We can bound the probability that $G$ contains $K_{s,t}$ by
\[
  \Pr[K_{s,t}\subset G]\leq \sum_{\substack{U\subset V\\\abs{U}=s}} \Pr\bigl[\abs{N(U)}\geq t\bigr]\leq 2\binom{n}{s}\frac{1}{t!}.
\]
(The two cases $U\subset L$ and $U\subset R$ are responsible for the factor of $2$.)
If $t= 10s \frac{\log n}{\log \log n}$, then this probability is very close to $0$.
Thus with high probability $G$ contains approximately $\Theta(n^{2-1/s})$ edges, and contains no $K_{s,t}$
for $t=10s \frac{\log n}{\log \log n}$.\medskip

The analysis above is nearly tight: one can show that for $t=0.1 s\frac{\log n}{\log \log n}$,
the random graph contains $K_{s,t}$ with overwhelming probability (see \cite[Section 4.5]{alon_spencer} for a proof of a 
similar result for cliques). 

So, the reason for the failure of probabilistic construction is that while $\abs{N(U)}$ has
mean $1$, the distribution of $\abs{N(U)}$ has a long, smoothly-decaying tail. Since there are many sets $U$,
it is likely that $\abs{N(U)}$ is large for some $U$.

\section*{Random algebraic construction}
Let $q$ be a prime power, and let $\F_q$ be the finite
field of order~$q$. We shall assume that $s\geq 4$ is fixed,
and that $q$ is sufficiently large as a function of $s$. 
Let $d=s^2-s+2$, $n=q^s$.  The graph $G$ that we will construct 
in this section will be bipartite. Each of the two parts, $L$ and $R$,
will be identified with $\F_q^s$.

Suppose $f$ is a polynomial in $2s$ variables over $\F_q$. We write the polynomial
as $f(X,Y)$ where $X=(X_1,\dots,X_s)$ and $Y=(Y_1,\dots,Y_s)$ are the first
and the last $s$ variables respectively. Such a polynomial induces a bipartite graph
in the natural way: pair $(x,y)\in L\times R$ is an edge if $f(x,y)=0$.

Let $\mathcal{P}\subset \mathbb{F}_q[X,Y]$ be the set of all polynomials of degree at most $d$ in each of $X$ and $Y$.
Pick a polynomial $f$ uniformly from $\mathcal{P}$ and let
$G$ be the associated graph. We shall show that $G$, on average, contains many
edges but hardly any copies of $K_{s,t}$ for $t=s^d+1$. We will then
remove few vertices from $G$ to render $G$ completely free of $K_{s,t}$'s while still leaving
many edges left.

We show that $G$ behaves very similarly to the random graph
that we constructed in the previous section with $p=1/q$. We begin by counting the number of edges in $G$.
\begin{lemma}\label{lem:exp}
For every $u,v\in\F_q^s$, we have $\Pr[f(u,v)=0]=1/q$.
In particular, the expected number of edges in $G$ is $n^2/q$.
\end{lemma}
\begin{proof}
Fix $u,v\in \F_q^s$. Let $\mathcal{P}_0=\{f\in\mathcal{P} : f(0,0)=0\}$ be the set of polynomials with zero constant term.
Every $f\in\mathcal{P}$ can be written uniquely as $f=g+h$, where $g\in \mathcal{P}_0$ 
and $h$ is a constant. So, a way to sample $f\in\mathcal{P}$ uniformly
is to first sample $g$ from $\mathcal{P}_0$, and then sample $h$ from $\F_q$. 
It is clear that having chosen $g$, out of $q$ possible choices for $h$ exactly one choice results in $f(u,v)=0$.
\end{proof}

To count the copies of $K_{s,t}$ we shall look at the distribution of $\abs{N(U)}$, where
$U$ is an arbitrary set of $s$ vertices in the same part. We shall focus on the case $U\subset L$, 
the other case being symmetric.

Computing the distribution of $\abs{N(U)}$ directly is hard. 
Instead we will compute moments of $\abs{N(U)}$ with aid of the following two lemmas:
\begin{lemma}\label{lem:onedim}
Suppose $u,u'\in \F_q^s$ are two distinct points, and $L$ 
is a linear function chosen uniformly among all linear functions $\F_q^s\to \F_q$.
Then $\Pr[Lu=Lu']=1/q$.
\end{lemma}
\begin{proof}
Since $u$ and $u'$ are distinct, there is a coordinate in which they differ. Without loss of generality,
it is the first coordinate. A linear function is uniquely determined by its action on
the basis vectors $e_1,\dotsc,e_s$. Sample $L$ by first sampling $Le_2,\dotsc,Le_s$ and then
sampling $Le_1$. Having chosen $Le_2,\dotsc,Le_s$ there is precisely one choice
for $Le_1$ such that $Lu=Lu'$.
\end{proof}
\begin{lemma}
Suppose $r,s\leq \min(\sqrt{q},d)$. Let $U\subset \F_q^s$ and $V\subset\F_q^s$ be sets of size $s$ and $r$ 
respectively. Then
\[
  \Pr\bigl[f(u,v)=0\text{ for all }u\in U,\, v\in V\bigr]=q^{-sr}.
\]
\end{lemma}
\begin{proof}
Call a set of points in $\F_q^s$ \emph{simple} if the first coordinates of all the points are distinct.

We first give the proof in the case when $U$ and $V$ are simple sets.
In this case, we decompose $f$ as $f=g+h$,
where $h$ contains the monomials $X_1^iY_1^j$ for $i=0,1,\dotsc,s-1$ and $j=0,1,\dotsc,r-1$,
whereas $g$ contains all the other monomials. Similarly to the proof of the preceding lemmas 
it is sufficient to show that the system of linear equations
\begin{equation}\label{eq:lineqs}
  h(u,v)=-g(u,v)\qquad\text{for all }u\in U,v\in V
\end{equation}
has a unique solution with polynomial $h$ as the unknown. This is a consequence of
the Lagrange interpolation theorem applied twice: the first application
yields, for each $u\in U$, single-variate polynomials $h_u(Y)$ of degree at most $r-1$
such that $h_u(v)=-g(u,v)$ for
all $v\in V$; the second application yields a polynomial $h(X,Y)$ such that each if the coefficients
of $h(u,Y)$ is equal to the respective coefficient of $h_u(Y)$ for all $u\in U$. That
latter condition implies of course that $h(u,v)=h_u(v)$. Note that the obtained polynomial $h$
is unique since the solution exists for each of $q^{rs}$ possible right-hand sides in \eqref{eq:lineqs},
and there are only $q^{rs}$ polynomials~$h$. 

We next treat the case of general $U$ and $V$. It suffices to find invertible linear
transformations $T$ and $S$ acting on $\F_q^s$ such that both $TU$ and $SV$ 
are simple. Indeed, the set of polynomials $\mathcal{P}$ is invariant
under change of coordinates in the first $s$ coordinates, and is invariant
under change of coordinates in the last $s$ coordinates. Hence, if we arrange
for $TU$ and $SV$ to be simple, we reduce to the special case treated above.

To find the requisite $T$, it suffices to find a linear map $T_1\colon \F_q^s\to\F_q$ that
is injective on~$U$. We can then find an invertible map $T\colon \F_q^s\to\F_q^s$ whose
first coordinate is $T_1$.
We pick $T_1$ uniformly at random from among all linear maps $\F_q^s\to\F_q$. By lemma \ref{lem:onedim}, 
for any distinct $x,x'\in X$, the probability that $T_1x=T_1x'$ is $1/q$, and so
\[
  \Pr\bigl[\exists x,x'\in X,\ x\neq x'\ \wedge\ T_1x=T_1x'\bigr]\leq \binom{s}{2}\frac{1}{q}<1,
\]
implying that a suitable $T_1$ (and hence $T$) exists. The construction of $S$ is analogous.
\end{proof}

Fix a set $U\subset \F_q^s$ of size $s$. For $v\in \F_q^s$, put $I(v)=1$
if $f(u,v)=0$ for all $u\in U$, and $I(v)=0$ if $f(u,v)\neq 0$ for some
$u\in U$.  The $d$'th moment of $\abs{N(U)}$ is easily computed by writing
$\abs{N(U)}$ as a sum of $I(v)$'s and expanding:
\begin{align*}
  \E\bigl[ \abs{N(U)}^d \bigr] &= \E\left[ \left(\sum_{v\in \F_q^s} I(v) \right)^d \right]
                               = \E\left[ \sum_{v_1,\dotsc,v_d\in \F_q^s} I(v_1)I(v_2)\dotsb I(v_d)\right]\\
                               &= \sum_{v_1,\dotsc,v_d\in \F_q^s} \E[I(v_1)I(v_2)\dotsb I(v_d)]\\
\intertext{%
The preceding lemma tells us that the summand is equal to $q^{-rs}$ if there are exactly $r$ distinct
points among $v_1,\dotsc,v_d$. Let $M_r$ be the number of surjective
functions from a $d$-element set onto an $r$-element set, and let $M=\sum_{r\leq d} M_r$. Breaking the sum according
to the number of distinct elements among $v_1,\dotsc,v_d$, we see that}
   \E\bigl[ \abs{N(U)}^d \bigr] &= \sum_{r\leq d} \binom{q^s}{r}M_r q^{-rs}
                               \leq \sum_{r\leq d} M_r =M.
\end{align*}
We can use the moments to bound the probability that $\abs{N(U)}$ is large:
\begin{equation}\label{eq:deviation}
  \Pr\bigl[\abs{N(U)}\geq \lambda\bigr]=\Pr\bigl[\abs{N(U)}^d\geq \lambda^d\bigr]\leq 
\frac{\E\bigl[ \abs{N(U)}^d \bigr]}{\lambda^d}\leq \frac{M}{\lambda^d}.
\end{equation}

We have shown that distribution of edges of $G$ enjoys some independence, and used that to derive \eqref{eq:deviation}.
It is now time to exploit the dependence between the edges of $G$.
The following result provides severe constraints on the values attainable by $\abs{N(U)}$:
\begin{lemma}\label{lem:bez}
For every $s$ and $d$ there exists a constant $C$ such the following holds:
Suppose $f_1(Y),\dotsc,f_s(Y)$ are $s$ polynomials on $\F_q^s$ of degree at most $d$, and consider
the set
\[
  W=\{ y\in \F_q^s : f_1(y)=\dotsb=f_s(y)=0 \}.
\]
Then exactly one of the following holds:
\begin{enumerate}
\item (Zero-dimensional case) $\abs{W}\leq C$,
\item (Higher-dimensional case) $\abs{W}\geq q-C\sqrt{q}$.
\end{enumerate}
The constant $C$ depends only on $s$ and the degrees of $f$'s.
\end{lemma}
\begin{proof}
The proof of this lemma is the sole place in the paper where we use algebraic geometry. A basic textbook is \cite{shafarevich}. 
For technical reasons, we will work not with projective, but with affine varieties, and so the intersection theory that we will 
employ differs slightly from the most common sources. Namely, we will use the results from \cite{heintz_affine_bezout}.
In particular, we use the same notion of the degree of a variety, namely $\deg V=\sum \deg V_i$, where the
sum is over irreducible components of $V$. The notion obeys the familiar properties:
First, the degree of the variety $\{f=0\}$, where $f$ is a non-zero polynomial, is at most $\deg f$. Second,
if $X$ is a zero-dimensional variety, then $\deg X$ is just the number of points in $X$.
Finally, in Theorem~1 on page 251 of the same paper, it is shown that the
Bezout's inequality holds, namely
\[
  \deg X\cap Y\leq \deg X\cdot \deg Y,
\]
for any two varieties $X,Y$. (A similar result in the projective space can
be found in \cite[Example~12.3.1]{fulton}.)

By dimension of a variety defined over a finite field $\F_q$, we will mean the dimension of
the variety as a variety over the algebraic closure $\overline{\F_q}$ (see \cite[Chapter 6]{shafarevich}).

To establish the lemma it suffices to prove that, for any fixed $m$ and $D$, whenever $V$ is an (affine) variety defined over $\F_q$ of degree $D$ and dimension $m$,
then the set of $\F_q$-points of $V$ satisfies either $\abs{V(F_q)}=O(1)$ or $\abs{V(F_q)}\geq q-O(\sqrt{q})$.
Indeed, $W$ is the set of $\F_q$-points of the variety with equations $f_1=\dotsb=f_s=0$,
and its degree is bounded by Bezout's inequality applied to the varieties $\{f_i=0\}$. Here ``bounded'' means
bounded in terms of $m$ and $D$; similarly, the constants in the big-oh notation are allowed to depend
on $m$ and~$D$.

The proof is by induction on $m$ (for all $D$ simultaneously). In the base case $m=0$ is trivial, as we then 
have $\abs{V}=D$. Suppose $m\geq 1$. If $V$ is reducible over $\F_q$, then
the degrees of the components add up to $D$, and we can treat each component
separately. So, assume that $V$ is irreducible over $\F_q$. If $V$ is also 
irreducible over $\overline{\F_q}$ then it has $q^{\dim V}(1-O(1/\sqrt{q}))$ points
by the Lang--Weil bound \cite{lang_weil} (for an elementary proof see \cite{schmidt}).
Otherwise $V$ is reducible over $\overline{\F_q}$, with $V_1,\dotsc,V_r$ as the components.
The reducibility means that $r\geq 2$.
The Frobenius automorphism $x\mapsto x^q$ acts on $V$, permuting the components. The action
is transitive because $V$ is irreducible over $\F_q$. Indeed, if $V_1,\dotsc,V_l$ is an orbit
of the action, then the variety $V_1\cup\dotsb\cup V_l$ is invariant under the action
of the Frobenius automorphism, and so is $\F_q$-definable \cite[Proof of Corollary~4]{tao}.
Similarly, each orbit gives rise to a proper subvariety of $V$.
As the union of these subvarieties is $V$, this contradicts the irreducibility of $V$.
The contradiction shows that the action is transitive, as claimed.

Let $V'=V_1\cap\dotsb\cap V_r$.
In view of the transitivity we have $V_1(\F_q)=\dotsb=V_r(\F_q)$, and so $V(\F_q)=V'(\F_q)$.
As $V'$ is invariant under the action of the Frobenius automorphism,
it is $\F_q$-definable. As $V$ is irreducible over $\F_q$, 
we cannot have $\dim V'=\dim V$ for it would follow that $V=V'$, contrary to $r\geq 2$..
Hence, $\dim V'<\dim V=m$. Moreover,
we can bound the degree of $V'$ via Bezout's inequality as follows
\begin{align*}
   \deg V'\leq \prod \deg V_i\leq (\tfrac{1}{r}\sum \deg V_i)^r=(D/r)^r\leq \exp(D/e).
\end{align*}
Since $V(\F_q)=V'(\F_q)$, the result follows from the induction hypothesis.
\end{proof}

We consider $s$ polynomials $f(u,\cdot)$ as $u$ ranges over $U$.
The preceding lemma then says that either $\abs{N(U)}\leq C$ or $\abs{N(U)}\geq q/2$
if $q$ is sufficiently large in terms of $s$. From \eqref{eq:deviation} we thus
obtain (for all sufficiently large $q$)
\[
  \Pr[\abs{N(U)}>C]=\Pr[\abs{N(U)}\geq q/2]\leq \frac{M}{(q/2)^d}.
\]

Call a set of $s$ vertices of $G$ \emph{bad} if their common neighborhood 
has more than $C$ vertices. Let $B$ the number of bad sets. The above
shows that 
\begin{equation}\label{eq:eb}
  \E[B]\leq 2\binom{n}{s}\frac{M}{(q/2)^d}=O(q^{s-2}).
\end{equation}

Remove a vertex from each bad set counted by $B$ from $G$ to obtain graph $G'$. 
Since no vertex has degree more than $q^s$, the number of edges in $G'$
is at most $Bq^s$ fewer than in $G$. Hence, the expected number of edges in $G'$
is at least
\[
  n^2/q-\E[B]q^s=\Omega(n^{2-1/s}),
\]
where $n^2/q$ comes from Lemma~\ref{lem:exp}, and the estimation of $\E[B]$ comes from \eqref{eq:eb}.

Therefore, there exists a graph with at most $2n$ vertices and $\Omega(n^{2-1/s})$ 
edges, but without $K_{s,C+1}$.\medskip

\begin{remark} 
An earlier version of this paper asserted that the constant $C$ in Lemma~\ref{lem:bez}
can be taken to be $\prod \deg f_i$. The assertion is false. Here is an example based on
the idea of Jacob Tsimerman. Let $a$ be any element of $\F_{p^2}$ that is not
in $\F_p$, and choose univariate polynomials $g$ and $h$ of degrees $d$ and $d-1$ respectively
that are completely reducible over $\F_p$ with distinct roots.
The bivariate polynomial $ag(x)+h(y)$ is irreducible over $\overline{\F_p}$.
Indeed, if it were reducible, then its Newton polygon\footnote{Newton polygon
of a bivariate polynomial $\sum_{i,j} a_{i,j} a_{i,j}x^i y^j$ is the convex
hull of $\{(i,j) : a_{i,j}\neq 0\}$.}
would be a Minkowski sum of Newton polygons of its factors \cite[Theorem VI]{ostrowski}.
Since $\operatorname\{(d,0),(0,d-1)\}$
is not a Minkowski sum of two smaller lattice polygons, $ag(x)+h(y)$ is irreducible over $\overline{\F_p}$.
Polynomial $a^p g(x)+h(y)$ is similarly
irreducible.  Let $f_1(x,y,z)=\bigl(ag(x)+h(y)\bigr)\bigl(a^p g(x)+h(y)\bigr)$. Since $f_1$
is invariant under the Frobenius automorphism, $f_1\in \F_p[x,y,z]$. Let $f_2(x,y,z)=f_3(x,y,z)=z$.
Then common zero set of $f_1,f_2,f_3$ is the set $\{(x,y,z) : g(x)=h(y)=z=0\}$ which has size $d(d-1)$, 
whereas $\prod \deg f_i=2d$.

\end{remark}

\textbf{Acknowledgements.} I am grateful to Roman Karasev for valuable discussions, and to David Conlon,
Zilin Jiang, and Eoin Patrick Long for comments on the earlier versions of this paper. I also thank the anonymous 
referee for detailed feedback and pointing reference \cite{heintz_affine_bezout}.

\bibliographystyle{plain}
\bibliography{turanconstr}

\begin{thebibliography}{10}

\bibitem{aigner_monthly}
Martin Aigner.
\newblock Tur\'an's graph theorem.
\newblock {\em Amer. Math. Monthly}, 102(9):808--816, 1995.

\bibitem{alon_ronyai_szabo}
Noga Alon, Lajos R{\'o}nyai, and Tibor Szab{\'o}.
\newblock Norm-graphs: variations and applications.
\newblock {\em J. Combin. Theory Ser. B}, 76(2):280--290, 1999.

\bibitem{alon_spencer}
Noga Alon and Joel~H. Spencer.
\newblock {\em The probabilistic method}.
\newblock Wiley-Interscience Series in Discrete Mathematics and Optimization.
  Wiley-Interscience [John Wiley \& Sons], New York, second edition, 2000.
\newblock With an appendix on the life and work of Paul Erd{\H{o}}s.

\bibitem{bbk}
Pavle V.~M. Blagojevi{\'c}, Boris Bukh, and Roman Karasev.
\newblock Tur\'an numbers for {$K_{s,t}$}-free graphs: topological obstructions
  and algebraic constructions.
\newblock {\em Israel J. Math.}, 197(1):199--214, 2013.

\bibitem{brown_constr}
W.~G. Brown.
\newblock On graphs that do not contain a {T}homsen graph.
\newblock {\em Canad. Math. Bull.}, 9:281--285, 1966.

\bibitem{erdos_renyi_sos}
P.~Erd{\H{o}}s, A.~R{\'e}nyi, and V.~T. S{\'o}s.
\newblock On a problem of graph theory.
\newblock {\em Studia Sci. Math. Hungar.}, 1:215--235, 1966.

\bibitem{erdos_simonovits}
P.~Erd{\H{o}}s and M.~Simonovits.
\newblock A limit theorem in graph theory.
\newblock {\em Studia Sci. Math. Hungar}, 1:51--57, 1966.

\bibitem{erdos_stone}
P.~Erd{\"o}s and A.~H. Stone.
\newblock On the structure of linear graphs.
\newblock {\em Bull. Amer. Math. Soc.}, 52:1087--1091, 1946.

\bibitem{fulton}
William Fulton.
\newblock {\em Intersection theory}, volume~2 of {\em Ergebnisse der Mathematik
  und ihrer Grenzgebiete. 3. Folge. A Series of Modern Surveys in Mathematics
  [Results in Mathematics and Related Areas. 3rd Series. A Series of Modern
  Surveys in Mathematics]}.
\newblock Springer-Verlag, Berlin, second edition, 1998.

\bibitem{furedi_bound}
Zolt{\'a}n F{\"u}redi.
\newblock An upper bound on {Z}arankiewicz' problem.
\newblock {\em Combin. Probab. Comput.}, 5(1):29--33, 1996.

\bibitem{heintz_affine_bezout}
Joos Heintz.
\newblock Definability and fast quantifier elimination in algebraically closed
  fields.
\newblock {\em Theoret. Comput. Sci.}, 24(3):239--277, 1983.
\newblock Also corrigendum in vol.~39 (1985), no. 2-3, 343.

\bibitem{kollar_ronyai_szabo}
J{\'a}nos Koll{\'a}r, Lajos R{\'o}nyai, and Tibor Szab{\'o}.
\newblock Norm-graphs and bipartite {T}ur\'an numbers.
\newblock {\em Combinatorica}, 16(3):399--406, 1996.

\bibitem{lang_weil}
Serge Lang and Andr{\'e} Weil.
\newblock Number of points of varieties in finite fields.
\newblock {\em Amer. J. Math.}, 76:819--827, 1954.

\bibitem{ostrowski}
A.~M. Ostrowski.
\newblock On multiplication and factorization of polynomials. {I}.
  {L}exicographic orderings and extreme aggregates of terms.
\newblock {\em Aequationes Math.}, 13(3):201--228, 1975.

\bibitem{schmidt}
Wolfgang~M. Schmidt.
\newblock {\em Equations over finite fields. {A}n elementary approach}.
\newblock Lecture Notes in Mathematics, Vol. 536. Springer-Verlag, Berlin-New
  York, 1976.

\bibitem{shafarevich}
I.~R. Shafarevich.
\newblock {\em Basic algebraic geometry}.
\newblock Springer-Verlag, New York-Heidelberg, 1974.
\newblock Translated from the Russian by K. A. Hirsch, Die Grundlehren der
  mathematischen Wissenschaften, Band 213.

\bibitem{simonovits_stability}
M.~Simonovits.
\newblock A method for solving extremal problems in graph theory, stability
  problems.
\newblock In {\em Theory of {G}raphs ({P}roc. {C}olloq., {T}ihany, 1966)},
  pages 279--319. Academic Press, New York, 1968.

\bibitem{tao}
Terence Tao.
\newblock Lang--{W}eil bound (blog post).
\newblock \url{http://terrytao.wordpress.com/2012/08/31/the-lang-weil-bound/}.

\bibitem{turan_orig}
P.~Tur{\'a}n.
\newblock On an extremal problem in graph theory (in {H}ungarian).
\newblock {\em Math. Fiz. Lapok}, 48:436--452, 1941.

\end{thebibliography}

\end{document}